\documentclass{article}

\usepackage{amssymb}
\usepackage{psfig}
\usepackage{amscd}

\usepackage[ruled]{algorithm}
\usepackage{algorithmic}

\newcommand{\rr}{{\mathbb R}}
\newcommand{\cc}{{\mathbb C}}
\newcommand{\ZZ}{{\mathbb Z}}

\newcommand{\bfc}{{\bf c}}

\newcommand{\bfh}{{\bf h}}

\newcommand{\bfl}{{\mbox{\boldmath $\lambda$}}}
\newcommand{\bfmu}{{\mbox{\boldmath $\mu$}}}

\newcommand{\x}{{\bf x}}
\newcommand{\y}{{\bf y}}
\newcommand{\zero}{{\bf 0}}

\newcommand{\p}{{\partial}}
\newcommand{\LC}{{\mathbf{lc}}}
\newcommand{\LE}{{\mathbf{le}}}
\newcommand{\LM}{{\mathbf{lm}}}
\newcommand{\LT}{{\mathbf{lt}}}
\newcommand{\IN}{{\mathbf{in}}}
\newcommand{\supp}{{\mbox{supp}}}
\newcommand{\corank}{{\mbox{corank}}}
\newcommand{\rank}{{\mbox{rank}}}

\newtheorem{theorem}{Theorem}[section]
\newtheorem{lemma}[theorem]{Lemma}
\newtheorem{proposition}[theorem]{Proposition}

\newcommand\qed{{\hspace*{\fill}$\Box$\vskip12pt plus 1pt}}

\begin{document}

\title{Newton's Method with Deflation for
       Isolated Singularities of Polynomial Systems}

\author{
Anton Leykin\thanks{
Department of Mathematics, Statistics, and Computer Science,
University of Illinois at Chicago, 851 South Morgan (M/C 249),
Chicago, IL 60607-7045, USA.
{\em Email:} leykin@math.uic.edu.
{\em URL:} http://www.math.uic.edu/{\~{}}leykin.}
\and
Jan Verschelde\thanks{
Department of Mathematics, Statistics, and Computer Science,
University of Illinois at Chicago, 851 South Morgan (M/C 249),
Chicago, IL 60607-7045, USA.
{\em Email:} jan@math.uic.edu or jan.verschelde@na-net.ornl.gov.
{\em URL:} http://www.math.uic.edu/{\~{}}jan.
This material is based upon work
supported by the National Science Foundation under Grant No.\
0105739 and Grant No.\ 0134611.}
\and
Ailing Zhao\thanks{
Department of Mathematics, Statistics, and Computer Science,
University of Illinois at Chicago, 851 South Morgan (M/C 249),
Chicago, IL 60607-7045, USA.
{\em Email:} azhao1@uic.edu.
{\em URL:} http://www.math.uic.edu/{\~{}}azhao1.}
}

\date{13 October 2004}

\maketitle

\begin{abstract}

\noindent  We present a modification of Newton's method
to restore quadratic convergence for isolated singular
solutions of polynomial systems.  Our method is symbolic-numeric:
we produce a new polynomial system which has the original
multiple solution as a regular root.
Using standard bases,
a tool for the symbolic computation of multiplicities,
we show that the number of deflation stages is bounded
by the multiplicity of the isolated root.
Our implementation performs well on a large class of
applications.



\end{abstract}

\section{Introduction}

Let $F(\x) = \zero$ be a polynomial system of $N$ equations
in $n$ unknowns $\x \in \cc^n$.
We are interested in $\x^*$, an {\em isolated} solution
of $F(\x) = \zero$:
\begin{equation}
   \mbox{\rm for small enough} ~ \epsilon > 0:
   \{ \ \y \in \cc^n: || \y - \x^* || < \epsilon \ \} \cap F^{-1}(\zero)
   = \{ \x^* \}.
\end{equation}
Denote by $A(\x)$ the Jacobian matrix of the system $F(\x) = \zero$.
We call $\x^*$ a {\em singular} solution
of $F(\x) = \zero$ $\Leftrightarrow$ $\rank(A(\x^*)) < n$.
Let $m$ be the {\em multiplicity} of the isolated solution $\x^*$
of $F(\x) = \zero$.

Newton's method (also called the method of Gauss-Newton when $N>n$)
generates a sequence of approximations $\x_k$ for $\x^*$.
If $\x^*$ is nonsingular, then the sequence converges quadratically
(i.e.: $\|\x_k-\x_{k+1} \| = O(\|\x_{k-1}-\x_k\|^2)$)
to~$\x^*$, which justifies its widespread usage.
But otherwise, if $\x^*$ is singular, the convergence slows down
and gets lost when~$\x_k \approx \x^*$.

The brutal force approach is to use a working precision of
$m \times D$ decimal places to achieve $D$ correct decimal
places in the final approximation.  Even as multiprecision
arithmetic is widely available and nowadays less expensive to use,
this approach can only work if all coefficients in the system $F$
have their first $m \times D$ decimal places correct.
Our goal is to restore the quadratic convergence of a sequence
converging to an isolated singular root without imposing
extra requirements of precision on $F$.

Newton's method for singular solutions has been extensively
researched.  The research up to the mid eighties is surveyed
in~\cite{Gri85}.  We classify research related to our work
in two domains:
\begin{description}
  \item[1. Detection and treatment of bifurcation points.]
When following a solution path of a system defined by a parameter,
the solution path may turn back or bifurcate for increasing values
of the parameter.  Techniques to detect and compute such bifurcation
points are generally done via Liapunov-Schmidt
reduction~\cite{AG03}~\cite{Gov00};
see also~\cite{LW93,LW94}, \cite{DFK87}, \cite{Kun89},
and the references cited therein.
  \item[2. Deflation method for polynomial systems.]
A symbolic deflation method was presented in~\cite{OWM83},
and further developed in~\cite{Oji87a}, \cite{Oji87b}, and~\cite{Oji88}.
We discovered this approach from the reference list of~\cite{Lec02}.
\end{description}
A theoretical framework to study the complexity and numerical
stability of Newton's method was developed by Shub and Smale,
see~\cite{BCSS98}, and was generalized to overdetermined systems
in~\cite{DS99}.  See~\cite{GLSY04b,GLSY04a} for recent generalizations
of this $\alpha$-theory to multiple roots.

The contribution of this paper is twofold.
First -- as announced in~\cite{VZ04} -- we provide a numerically
stable implementation of a modified symbolic deflation method.
Second, using standard bases~\cite{GP96}, we show that the number
of deflations needed to restore the quadratic convergence of
Newton's method is bounded by the multiplicity.
In the next section we describe our method,
followed by an introduction to standard bases and our proof
in the third section.  Our symbolic-numeric implementations and
numerical results are described in sections four and five.

\noindent {\bf Acknowledgements.}  Our modified method with numerical
results was presented at a poster session at ISSAC'04, 6 July 2004,
Santander, Spain.  We are grateful for the reactions of Gr\'egoire
Lecerf -- in particular for pointing at the work of Kunkel -- and
to Erich Kaltofen and Zhonggang Zeng.  Also discussions with Andrew
Sommese and Charles Wampler at Oberwolfach helped us.

\section{A Modified Deflation Method}

A singular root $\x^*$ of a square (i.e.: $N=n$)
system $F(\x) = \zero$ satisfies
\begin{equation} \label{eq_augmented}
\left\{
     \begin{array} {c}
        F(\x) = \zero \\
        \det(A(\x)) = 0.
     \end{array}
  \right.
\end{equation}
The augmented system~(\ref{eq_augmented})
forms the basic idea for deflation.
If $\x^*$ is isolated and $\corank(A(\x^*)) = 1$,
then $\x^*$ as root of~(\ref{eq_augmented}) has
a lower multiplicity.

We find deflation used repeatedly first in~\cite{OWM83},
and later modified in~\cite{Oji87a} and applied in~\cite{Oji87b,Oji88}.

In theory, $\det(A(\x))=0$ (or maximal minors)
could be used to form new equations.
But this is neither good symbolically because the determinant is usually
of high degree and leads to expression swell, nor numerically, as the
evaluation of the determinant is numerically unstable:
$\| \det(A(\x)) - \det(\overline{A}(\x)) \|
\gg \| A(\x) - \overline{A}(\x) \|$.

Instead of using the determinant,
on a system $F$ of $N$ equations in $n$ variables,
we proceed along the
following three steps to form new equations:
\begin{enumerate}
\item Let $r = \rank(A(\x_0))$ for $\x_0 \approx \x^*$.
      For numerical stability, we compute the rank via a Singular Value
      Decomposition (SVD) of the matrix~$A = A(\x_0)$.
\item Let $\bfh \in \cc^{r+1}$ be a random vector
      and $B \in \cc^{n\times(r+1)}$ be a random matrix.
      For numerical stability, we generate random numbers
      on the unit circle.
\item Let $C(\x) = A(\x) B$, notice that $C = [\bfc_1, \mathbf
c_2, \ldots, \mathbf c_{r+1}]$ is an $N\times(r+1)$ matrix with
polynomial entries.

\end{enumerate}
With probability one (exceptional pairs of vectors $\bfh$ and
matrices $B$ belong to a proper algebraic subset of
$\cc^{r+1}\times\cc^{n\times(r+1)}$) we have
\begin{eqnarray*}
  \rank(A(\x^*)) = r
  & \Leftrightarrow &
  \corank(C(\x^*)) = 1\\
  & \Leftrightarrow & \exists\,!\ \bfl =
      \left(\begin{array}{c}\lambda_1\\\lambda_2\\\vdots\\\lambda_{r+1}\end{array}\right):
G(\x^*, \bfl)=
   \left(
     \begin{array}{c}
       \displaystyle \sum_{i=1}^{r+1} \lambda_i \bfc_i(\x^*) \\
       \displaystyle \sum_{i=1}^{r+1} h_i \lambda_i  - 1
     \end{array}
   \right) = \zero.
\end{eqnarray*}

The random $\bfh$ and $B$ guarantee the existence and uniqueness
of the solution $\bfl$ to $G(\x, \bfl)$ when $\x=\x^*$.

In one deflation step, we add the equations
of $G(\x, \bfl)$ instead of $\det(A(\x))=0$ to the
system $F(\x) = \zero$, adding $r+1$ extra variables $\lambda_1$,
$\lambda_2$, $\ldots$, $\lambda_{r+1}$.

For $r = n-1$, adding $G(\x,\bfl)$ to $F(\x) = \zero$ is
equivalent to the Liapunov-Schmidt reduction~\cite{AG03} used
to compute bifurcation points, see e.g.~\cite{DFK87}.

\begin{figure}[hbt]
\begin{center}
\begin{picture}(300,310)(30,0)

\put(175,285){\oval(300,40)}
\put( 75,290){Input: $F(\x) = \zero$ polynomial system;}
\put(105,275){$\x_0$ initial approximation for $\x^*$.}

\put(175,265){\vector(0,-1){15}}

\put(80,210){\framebox(190,40)[l]}

\put(95,235){$A^+ := {\rm SVD}(A(\mathbf x_k))$;}

\put(95,220){$\x_{k+1} := \x_k - A^+ F(\x_k)$;}

\put(175,210){\vector(0,-1){15}}

\put(175,195){\line(2,-1){60}}

\put(175,195){\line(-2,-1){60}}

\put(175,135){\line(-2,1){60}}

\put(175,135){\line(2,1){60}}

\put(154,170){Quadratic}
\put(150,155){convergence?}

\put(235,172){Yes}

\put(235,165){\vector(1,0){22}}

\put(305,165){\oval(95,40)}

\put(265,162){Output: $F; \x_{k+1}$.}

\put(175,135){\vector(0,-1){15}}

\put(180,125){No}

\put(55,15){\framebox(240,105)[l]}

\put(70,100){$r:= \rank(A^+)$;}

\put(70,73){$F:=\left\{\begin{array} {c}
      F(\x) = \zero\\
      G(\x, \bfl) = \zero
   \end{array}\right.$;}

\put(70,45){$\bfl := {\rm LeastSquares} (G(\x_{k+1},\bfl))$;}

\put(70,30){$k:=k+1$; $\quad \x_k:=(\x_k, \bfl)$;}

\put(175,15){\line(0,-1){15}}

\put(25,0){\line(1,0){150}}

\put(25,0){\line(0,1){230}}

\put(25,230){\vector(1,0){55}}

\end{picture}
\caption{Flowchart for a modified deflation method.
In practice, the test for quadratic convergence is implemented
by a test whether the numerical rank of $A(\x_{k})$ equals
the number of columns of~$A(\x_{k})$.}
\label{figflowchart}
\end{center}
\end{figure}
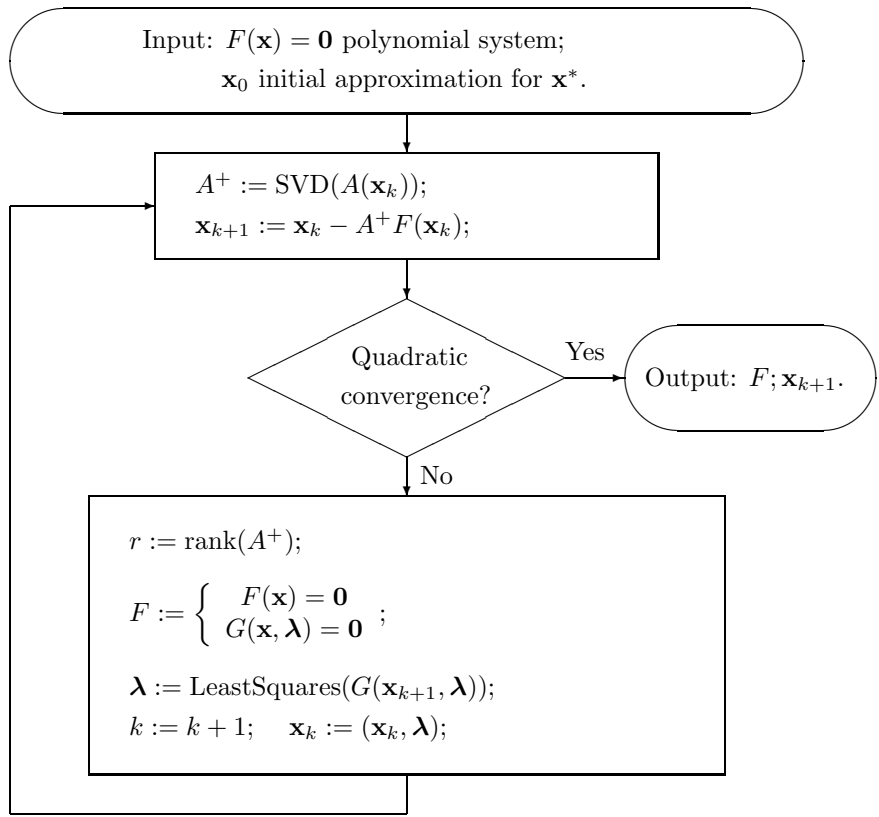

\newpage
\section{A Bound on the Number of Deflations}

The termination of our algorithm in Figure~\ref{figflowchart}
depends on the following theorem.

\begin{theorem}
The number of deflations needed to restore the quadratic convergence
of Newton's method converging to an isolated solution is strictly
less than the multiplicity of the isolated solution.
\end{theorem}

The answer to the question ``How much less?'' can be understood
by looking at a standard basis for the ideal generated by the
given polynomials in the system.  We use standard bases to prove
the termination of our algorithm, as explained in the next two
subsections.

\subsection{Standard Bases for Local Orderings}

Let $R=k[x_1,...,x_n]$ be the ring of polynomials in $n$ variables
with coefficients in the field~$k$.
We use the following multi-degree notation:
$\x^\alpha = x_1^{\alpha_1} x_2^{\alpha_2} \cdots x_n^{\alpha_n}$,
where $\alpha=(\alpha_1,...,\alpha_n)$ is a
vector of nonnegative integers.

A multiplicative ordering $\leq$ on the monoid $\{ \x^\alpha\ |\
\alpha\in \ZZ_{\geq 0}^n \}$ is a \emph{local ordering} if $1 >
\x^\alpha$ for all $\alpha \neq (0,0,\ldots,0)$.

To any \emph{weight vector} $\omega \in \rr_{< 0}^n$ we may
associate the weight ordering $\leq_\omega$ by setting
\begin{equation}
  \x^\alpha \leq_\omega \x^\beta
  \Longleftrightarrow \langle \alpha , \omega \rangle \leq
                  \langle \beta , \omega \rangle\mbox{ \textbf{or} } [\ \langle \alpha , \omega \rangle
                  =
                  \langle \beta , \omega \rangle \mbox{ \textbf{and} } \x^\alpha\leq_{lex}\x^\beta\ ],
\end{equation}
where $\langle \cdot , \cdot \rangle$ is the usual inner product
and the lexicographic ordering, $\leq{lex}$, is used to break
ties.

In presence of a monomial ordering $\leq_\omega$, a polynomial
\begin{equation}
  f(\x) = \sum_{\alpha \in \ZZ_{\geq 0}^n} c_\alpha \x^\alpha \in R
  \mbox{, where $\supp(f) = \{ \ \alpha \ | \ c_\alpha\neq 0 \ \}$ is finite,}
\end{equation}
has the following attributes associated with it:
\begin{eqnarray*}
  \LE(f) &=& \mbox{the leading exponent} = \max_{\leq_\omega} \supp(f) \\
  \LM(f) &=& \mbox{the leading monomial} = \x^{\LE(f)} \\
  \LC(f) &=& \mbox{the leading coefficient} = c_{\LE(f)} \\
  \LT(f) &=& \mbox{the leading term} = \LC(f)\LM(f)
\end{eqnarray*}

Let $I \subset R$ be an ideal.  We call a set of polynomials
$S \subset I$ a \emph{standard basis} of $I$ if for any $f\in I$
there is $g \in S$ such that $\LM(g)|\LM(f)$. Alternatively, $S$
is a standard basis iff the \emph{initial ideal}
$\IN(I) = \langle \{\LM(f) | f \in I\} \rangle$
is generated by the leading monomials
$\LM(S) = \{ \LM(g) | g \in S \}$.

The monomials that do not belong to the initial ideal $\IN(I)$ are
called \emph{standard monomials}.  The minimal generators of
$\IN(I)$ shall be called the \emph{corners} of the staircase, the
corners of the form $x_i^a$ for some~$i$ and~$a$ are called the
\emph{endpoints} of the staircase.

A standard basis $S$ is called \emph{reduced} if the leading
monomials of its elements form a minimal generating set for the
initial ideal $\IN(I)$ and the tail $g-\LT(g)$ contains only
standard monomials.

Graphically, any monomial ideal can be represented by a staircase
in the nonnegative integer lattice $\ZZ^n_{\geq 0}$. For example,
let $I$ be the ideal of $R=k[x_1,x_2]$ generated by
\begin{eqnarray}
f_1&=&x_1^3+x_1 x_2^2; \nonumber\\
f_2&=& x_1 x_2^2 + x_2^3; \\
f_3&=& x_1^2x_2+x_1x_2^2. \nonumber \label{exa1}
\end{eqnarray}

\begin{figure}[hbt]
\unitlength 0.9mm
\begin{picture}(160,50)(0,0)

\put(5,0){
\begin{picture}(80,50)(0,0)

\linethickness{0.3mm} \put(0,0){\line(0,1){50}}
                      \put(0,50){\vector(0,1){0.12}}
\linethickness{0.3mm} \put(0,0){\line(1,0){50}}
                      \put(50,0){\vector(1,0){0.12}}
\linethickness{0.3mm} \put(0,0){\circle{4}}
\linethickness{0.3mm} \put(0,10){\circle{4}}
\linethickness{0.3mm} \put(20,0){\circle{4}}
\linethickness{0.3mm} \put(0,20){\circle{4}}
\linethickness{0.3mm} \put(10,10){\circle{4}}
\linethickness{0.3mm} \put(10,0){\circle{4}}
\linethickness{0.3mm} \put(0,40){\circle*{4}}
\linethickness{0.3mm} \put(20,10){\circle*{4}}
\linethickness{0.3mm} \put(30,0){\circle*{4}}
\linethickness{0.6mm} \put(0,40){\line(1,0){10}}
\linethickness{0.6mm} \put(20,10){\line(0,1){10}}
\linethickness{0.6mm} \put(20,10){\line(1,0){10}}
\linethickness{0.6mm} \put(30,0){\line(0,1){10}}

\put(14,24){\makebox(0,0)[cl]{$\underline{x_1x_2^2}+x_2^3$}}
\put(66,14){\makebox(0,0)[cc]{}}
\put(24,14){\makebox(0,0)[cl]{$\underline{x_1^2x_2}+x_1x_2^2$}}
\put(34,4){\makebox(0,0)[cl]{$\underline{x_1^3}+x_1x_2^2$}}

\linethickness{0.6mm} \put(10,20){\line(0,1){20}}
\linethickness{0.6mm} \put(10,20){\line(1,0){10}}
\linethickness{0.3mm} \put(10,20){\circle*{4}}
\linethickness{0.3mm} \put(0,30){\circle{4}}

\put(2,44){\makebox(0,0)[cl]{$\underline{x_2^4}$}}
\put(8,44){\makebox(0,0)[cl]{}}
\put(24,34){\makebox(0,0)[cl]{$\omega$}}
\put(0,20){\makebox(0,0)[cl]{}}


\put(22,32){\makebox(0,0)[cl]{}}

\linethickness{0.3mm}
\multiput(20,30)(0.12,0.24){83}{\line(0,1){0.24}}
\put(20,30){\vector(-1,-2){0.12}}
\end{picture}
}

\put(75,0){
\begin{picture}(80,50)(0,0)
\linethickness{0.3mm}
\put(0,0){\line(0,1){50}}
\put(0,50){\vector(0,1){0.12}}
\linethickness{0.3mm}
\put(0,0){\line(1,0){50}} \put(50,0){\vector(1,0){0.12}}
\linethickness{0.3mm} \put(0,0){\circle{4}}
\linethickness{0.3mm} \put(0,10){\circle{4}}
\linethickness{0.3mm} \put(20,0){\circle{4}}
\linethickness{0.3mm} \put(0,20){\circle{4}}
\linethickness{0.3mm} \put(10,10){\circle{4}}
\linethickness{0.3mm} \put(10,0){\circle{4}}
\linethickness{0.3mm} \put(0,30){\circle*{4}}
\linethickness{0.3mm} \put(20,10){\circle*{4}}
\linethickness{0.3mm} \put(40,0){\circle*{4}}
\linethickness{0.6mm} \put(0,30){\line(1,0){10}}
\linethickness{0.6mm} \put(20,10){\line(0,1){10}}
\linethickness{0.6mm} \put(20,10){\line(1,0){20}}
\linethickness{0.6mm} \put(40,0){\line(0,1){10}}

\put(24,14){\makebox(0,0)[cl]{$\underline{x_1x_2^2}-x_1^3$}}
\put(66,14){\makebox(0,0)[cc]{}}
\put(14,24){\makebox(0,0)[cl]{$\underline{x_1^2x_2}+x_1x_2^2$}}
\put(4,34){\makebox(0,0)[cl]{$\underline{x_2^3}+x_1x_2^2$}}

\linethickness{0.6mm} \put(10,20){\line(0,1){10}}
\linethickness{0.6mm} \put(10,20){\line(1,0){10}}
\linethickness{0.3mm} \put(10,20){\circle*{4}}
\linethickness{0.3mm} \put(30,0){\circle{4}}

\put(44,4){\makebox(0,0)[cl]{$\underline{x_1^4}$}}
\put(8,44){\makebox(0,0)[cl]{}}
\put(28,32){\makebox(0,0)[cl]{$\omega$}}
\put(0,20){\makebox(0,0)[cl]{}}


\put(22,32){\makebox(0,0)[cl]{}}

\linethickness{0.3mm}
\multiput(20,30)(0.24,0.12){83}{\line(1,0){0.24}}
\put(20,30){\vector(-2,-1){0.12}}
\end{picture}
}
\end{picture}
\caption{Two different staircases of the standard basis of $I$
         with respect to different local orderings
         $\leq_{(-1,-2)}$ (at the left)
         and $\leq_{(-2,-1)}$ (at the right).
Monomials generating $\IN(I)$ are represented by black disks,
while the standard monomials are shown as empty circles.}
\label{figstaircases}
\end{figure}
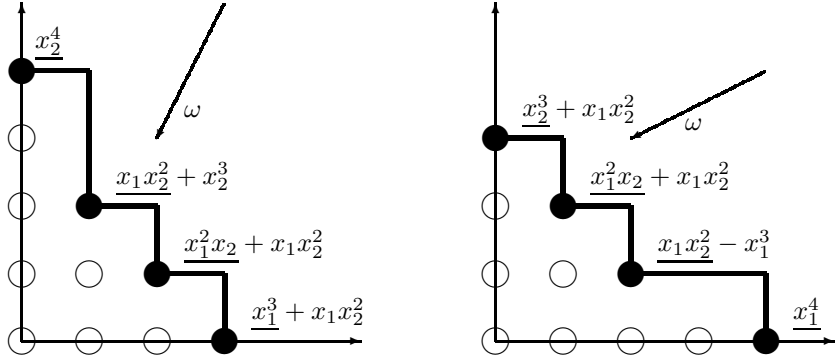

The initial ideal depends on the ordering chosen:  the staircase
at the left in Figure~\ref{figstaircases} represents $\IN_w(I)$,
where $\omega = (-1,-2)$. The staircase at the right in
Figure~\ref{figstaircases} represents $\IN_w(I)$, where $\omega =
(-2,-1)$.

Observe in Figure~\ref{figstaircases} that the number of standard
monomials is the same for both orderings. This is always so, for
any local ordering, as the standard monomials form a basis of the
$k$-linear space $R_{\langle x_1, \ldots,x_n \rangle}/R_{\langle
x_1, \ldots ,x_n \rangle}I$, where $R_{\langle x_1, \ldots ,x_n
\rangle}$ is the localization of the polynomial ring $R$ at the
origin and $R_{\langle x_1, \ldots ,x_n \rangle}I$ is the
extension of the ideal $I$ in this localized ring
(see~\cite{CLO98} for details). This linear space is of finite
dimension iff the origin is an isolated solution; its dimension,
which is the multiplicity of the origin, then equals the number of
the standard monomials for any local ordering.

Thanks to Mora~\cite{Mor82}, there is an algorithm for computing
standard bases, generalized by Greuel and Pfister~\cite{GP96}, and
implemented in the computer algebra system Singular~\cite{GP02}.
We will not use standard bases explicitly -- except for
theoretical purposes -- but note an analytic interpretation of the
local ordering $\leq_w$: as we approach the origin along a smooth
curve
\begin{equation}
  c: \cc \rightarrow \cc^n \quad {\rm such~that} \quad
  c(t)=
  \left\{
    \begin{array}{c}
       b_1 t^{-\omega_1}(1+O(t)) \\ \vdots \\
       b_n t^{-\omega_n}(1+O(t))
    \end{array}
  \right.
\end{equation}
with $\omega \in \ZZ_{< 0}^n$ and $(b_1,...,b_n) \in \cc^n
\setminus \{ 0 \}$, for every $f\in I$ the leading term $\LT_w(f)$
becomes dominant, i.e.:
\begin{equation}
f(c(t))=\LT_\omega(f)(c(t)) + O(t^{-\langle \omega, \LE(f)\rangle}).
\end{equation}

\subsection{Understanding the Deflation Method}

First of all let us formulate the goal of what we would call the
\emph{symbolic deflation} process: Given a system of polynomial
equations $f_i=0$, $i=1,2,\ldots, N$ with the point $\x^* \in \cc^n$
as an isolated solution of multiplicity $m>1$, find a system
$g_i=0$, $i=1,2, \ldots ,N'$ such that $\x^*$ is still an isolated
solution of multiplicity less than~$m$.

The best deflation one can cook up is the one that corresponds to
the maximal ideal annihilating $x^*_i$, i.e.: $g_i=x_i-x^*_i$,
$i=1,2,\ldots,n$. However, from a practical angle of numerical
methods what we actually need is an algorithm that would relate
the deflated system to the original one in a numerically stable
way and taking into account the fact that the isolated solution
$\x^*$ may be known only approximately.

\subsubsection{A Symbolic Deflation Method}

Here we assume that everything is exact and, therefore, w.l.o.g.
we may assume that the isolated solution $\x^*$ is the origin.

Consider the ideal $I$ generated by the polynomials $f_i$ of the
original system. We call an ideal $I'$ a \emph{deflation} of $I$
if $I'\supset I$, $I'\neq R$, and the multiplicity of the origin
for $I'$ is lower than that for the original ideal $I$.

If the multiplicity $m>1$, it means that the initial ideal
$\IN(I)$ does not contain $x_i$ for some~$i$.

\begin{proposition}\label{propDeflate}
Suppose $m > 1$ and let $g$ be an element of a reduced standard
basis of $I$ with respect to a local monomial ordering $\leq$,
such that $\LM(g) = x_i^d$, for some $i\in\{1,...,n\}$ and $d >
1$. Then the ideal $I' = I+ \langle \p g/\p x_i \rangle$ is a
deflation of~$I$.
\end{proposition}

\noindent {\em Proof.} The derivative $\p g/\p x_i$ can not
contain monomials $> x_i^{d-1}$. Therefore, $I'$ contains $I$
properly, since $\LM(\p g/\p x_i)=x_i^{d-1}$ is a standard monomial
for~$I$. The appended generator $\p g/ \p x_i$ still vanishes at
the origin, hence, $I'\neq R$.~\qed

Life is easy if one can compute a standard basis. However, this
procedure is expensive symbolically and unstable numerically. Can
we find $x_i$ and $g$ in the proposition in a less straightforward
way? The next lemma gives a positive answer.

A linear coordinate change $T:\cc^n\rightarrow \cc^n$ induces an
automorphism of the polynomial ring $R=\cc[x_1,...,x_n]$, which we
call $T$ as well: $T(f)(x) = f(T(x))$. The ideal
$T(I)=\{T(f)~|~f\in I \}=\langle T(f_1),\ldots,T(f_N) \rangle$
represents the system after the change of coordinates.

Let $A(\x)$ be the Jacobian matrix of the system $F(\x)=0$, i.e.:
an $N$-by-$n$ matrix with polynomial entries $A_{ij}(\x) = \p f_i
/ \p x_j$. The origin is singular iff $c = \corank(A(\zero))>0$.
Since the Jacobian matrix is rank deficient, the kernel of
$A(\zero)$ is nonzero.

\begin{lemma} \label{lemmaCoordChange}
Take a nonzero vector $\bfl \in \ker A(\zero) \subset \cc^n$ and
let $T:\cc^n\rightarrow \cc^n$ be a linear coordinate
transformation such that:
\begin{equation}
  T_i(\x) = \lambda_i x_1 + \sum_{j=2}^{n} \mu_{ij} x_j, \quad
  \mbox{~for~} i = 1,2,\ldots,n,
\end{equation}
where $[\bfl,\bfmu_2,...,\bfmu_n]$ is a nonsingular matrix.

Then $\p_1 (T(I)) = \{\frac{\p}{\p x_1} f~|~f\in T(I)\}$ is a
deflation of $T(I)$.
\end{lemma}

\noindent {\em Proof.} For all $i=1,2,\ldots,N$,
\begin{equation}\label{equDiffWRT1}
\frac{\p}{\p x_1} (f_i(T(\x)))\\
  = \sum_{j=1}^n \frac{\p f_i}{\p x_j} (T(\x))
           \cdot \frac{\p T_j}{\p x_1} (\x)
  = \sum_{j=1}^n \left( \frac{\p f_i}{\p x_j} (T(\x)) \right)
  \lambda_j.
\end{equation}
The last expression is equal to $0$ when $\x=\zero$, since $\bfl
\in \ker A(\zero)$.

Take any $g = b_1 T(f_1)+ \cdots + b_N T(f_N) \in T(I)$, where
$b_i \in R$ for all $i$. Then
\begin{eqnarray*}
  \frac{\p g}{\p x_1}
           & = & b_1 \frac{\p (T(f_1))}{\p x_1}  + \cdots
               + b_N \frac{\p (T(f_N))}{\p x_1}  \\
           & + & T(f_1) \frac{\p b_1}{\p x_1}  + \cdots
               + T(f_N) \frac{\p b_N}{\p x_1}.
\end{eqnarray*}
In view of (\ref{equDiffWRT1}), the last expression evaluates to
$0$ at $\x=\zero$. Therefore, $\p_1(T(I))$ is a proper ideal
annihilating the origin. On the other hand, there is an element
$g$ of a reduced standard basis of $T(I)$ with respect to a local
ordering such that $\IN(g)=x_1^d$ with $d>1$. According to the
Proposition \ref{propDeflate} the ideal $I' = T(I)+ \langle \p
g/\p x_1 \rangle$ is a deflation of $T(I)$. So is $\p_1(T(I))$,
for it contains $I'$. ~\qed

Lemma~\ref{lemmaCoordChange} leads
to Algorithm~\ref{Alg_Symbolic_Deflation}.

\begin{algorithm}[hbt]
\caption{G = Symbolic\_Deflation($F$)}
\label{Alg_Symbolic_Deflation}

\begin{algorithmic}
\REQUIRE $F$, a finite set of polynomials in $R$, such that the
         ideal $\langle F \rangle$ has multiplicity $m>1$ at the origin.

\ENSURE $G$, a finite set of polynomials in $R$, such that the
        ideal $\langle G \rangle + \langle F \rangle$
        is a deflation of $\langle F \rangle$.

\smallskip \hrule \smallskip

\STATE Compute the Jacobian $A$ of $F$ at the origin;

\STATE Pick a nonzero vector $\bfl \in \ker A(\zero)$;

\STATE $\begin{displaystyle}
          G := \left\{
              \sum_{i=1}^{n} \lambda_i \frac{\p f}{\p x_i}
                      \ | \ f\in F
               \right\}
        \end{displaystyle}$.

\end{algorithmic}
\end{algorithm}

\subsubsection{A Numerical Deflation Method}

\begin{algorithm}[hbt]
\caption{G = Numeric\_Deflation($F,\x_0$)}
\label{Alg_Numeric_Deflation}

\begin{algorithmic}
\REQUIRE $F=\{f_1,\ldots,f_N\}$, a finite set of polynomials in $R$,
such that the ideal $\langle F \rangle$ has multiplicity $m>1$
at the point $\x^* \approx \x_0$.

\ENSURE $G$, a finite set of polynomials in
   $R'=R[\lambda_1,\ldots,\lambda_{r+1}]$, where $r = \rank A(\x_0)$,
   such that \\
$\bullet$ the ideal $\langle G \rangle \subset R'$ has an isolated
solution at
the point $P=(\x^*,\bfl^*) \in \cc^{n,r+1}$;\\
$\bullet$ the vector $\bfl^*$ is determined uniquely; \\
$\bullet$ the multiplicity of $P$ is less than $m$.

\smallskip \hrule \smallskip

\STATE Compute the Jacobian matrix $A(\x)$ of $F$;

\STATE $r := \rank A(\x_0)$;

\STATE \textbf{(R1)} Generate a random matrix
                     $B \in \cc^{n \times (r+1)}$;

\STATE $C(\x):= A(\x) B$
($N\times (r+1)$ matrix with polynomial entries);

\STATE Let $\bfl = (\lambda_1,\ldots,\lambda_{r+1})^T$
       be a vector of indeterminates;

\STATE Consider $N$ new polynomials
       $g_i(\x,\bfl) = (C(\x)\bfl)_i \in R'$;

\STATE \textbf{(R2)}
    $h(\bfl) :=
     h_1 \lambda_1 + \cdots + h_{r+1} \lambda_{r+1} - 1$,
where the $h_i$ are random numbers in $\cc$;

\STATE $G := F \cup \{ g_1, \ldots, g_N \} \cup \{ \bfh \}$.
\end{algorithmic}
\end{algorithm}

This method (see Algorithm~\ref{Alg_Numeric_Deflation})
does boil down to the
Algorithm~\ref{Alg_Symbolic_Deflation} and here is why.

Consider a point $P = (\x,\bfl) \in \cc^{n,r+1}$
and let $\x=\x_0$.
When this specialization is performed, the values for $\bfl$
are determined by the following system of $N+1$ linear equations:
\begin{equation}
  \left\{
     \begin{array}{rll}
        A(\x_0) \bfl & = & 0\\
       \langle \bfh , \bfl \rangle & = & 1
     \end{array}
  \right.
\end{equation}
where $\bfh = (h_1, \ldots, h_{r+1})$
is a vector of random complex numbers.

Observe how $C$ is created: the randomization step \textbf{R1}
insures that the $r+1$ columns of $C(\x^*)$ are random combinations
of the columns of $A(\x^*)$ and, therefore, $\corank C(\x^*) = 1$
with probability one. Then $\bfl$ is bound to live in the one-dimensional
$\ker C(\x^*)$. The randomization step \textbf{R2} makes sure one
nonzero vector is picked out from the kernel. This proves the
uniqueness of $\bfl^*$.

Since $C(\x^*) \bfl = A(\x^*) B \bfl = 0$,
the vector $B \bfl$ corresponds to $\bfl$
in Algorithm~\ref{Alg_Symbolic_Deflation},
provided $\x^*$ is the origin.
Therefore, the multiplicity drops by Lemma~\ref{lemmaCoordChange}.

\section{A Symbolic-Numeric Implementation}

The method was tested and developed in Maple~9.
It is implemented in PHCpack~\cite{V99},
publicly available in release 2.3.
While floating-point arithmetic is used, the result is symbolic,
in the form of a new polynomial system with a well conditioned root.

In our implementation we computed the full SVD, although
we only need to know the rank of the Jacobian matrix at
the current approximation for the root.  To speed up the
algorithm, we refer to the recent techniques presented in~\cite{LZ04}.

Another issue is the exploitation of the structure
in the deflated systems.  This issue is related to
automatic differentiation~\cite{Gri00,Ral81}.
See also~\cite{Kun96}.

The deflation of the system $F(\x) = \zero$
with Jacobian matrix $A(\x)$ leads to a system
$F^{(1)}(\x) = \zero$ whose Jacobian matrix
$A^{(1)}(\x,\bfl)$ has an obvious block structure.
At the right of formula~(\ref{eqblock1}), the matrix
$A^{(1)}(\x,\bfl)$ has two columns: the first with
derivatives with respect to~$\x$ and the second with
derivatives with respect to~$\bfl$.

\begin{equation} \label{eqblock1}
     F^{(1)}(\x,\bfl)
     \! = \! \left\{
          \begin{array}{lcr}
             F(\x) & \!\! = \!\! & \zero \\
             A(\x) B \bfl & \!\! = \!\!  & \zero \\
             \multicolumn{1}{r}{\bfh \bfl} & \!\! = \!\!  & 1
          \end{array}
       \right.
  \quad
     A^{(1)}(\x,\bfl)
     \! = \! \left[
       \begin{array}{cc}
          A(\x) & \zero \\
          \left( \frac{\partial}{\partial \x} A(\x) \right) B \bfl &
          A(\x) B \\
          \zero & \bfh
       \end{array}
    \right]
\end{equation}
where the derivative of a matrix of polynomials $A(\x)$
with respect to a vector of unknowns $\x = (x_1,x_2,\ldots,x_n)$
gives a sequence of matrices:
\begin{equation} \label{eq_derivative_of_matrix}
  \frac{\partial}{\partial \x} A(\x)
  = \left[
          \frac{\partial A(\x)}{\partial x_1} \quad
          \frac{\partial A(\x)}{\partial x_2} \quad \cdots \quad
          \frac{\partial A(\x)}{\partial x_n}
    \right].
\end{equation}
To evaluate $A^{(1)}(\x,\bfl)$ efficiently, we must first
evaluate $A(\x)$ {\em before} we multiply with $B$.
Otherwise, multiplying first the polynomials in $A(\x)$ with
the matrix~$B$ will give rise to multiple occurrences of the
same polynomials in~$A(\x)$.  Expanding $A(\x) B$ as polynomials
in~$\x$ makes matters even worse.
Likewise, we first evaluate $\frac{\partial}{\partial \x} A(\x)$
{\em before} multiplying with the vector~$B \bfl$.

With the exploitation of structure, the system and its Jacobian
matrix can be evaluated about twice as fast, as our experiences
on the cyclic 9-roots problem suggest.  With 7 multipliers added
to the original 9 variables, the time for 1000 evaluations
of $F^{(1)}(\x,\bfl)$ at a
random point drops from 1.75 to 0.94 cpu seconds\footnote{Execution
done on a 2.4GHz Linux machine with full optimization of the code.},
and the time to evaluate $A^{(1)}(\x,\bfl)$ 1000 times drops
from 11.0 to 5.3 cpu seconds.
With 8 multipliers, the drop is even more significant:
from 2.07 to 0.92 cpu seconds for 1000 evaluations of~$F^{(1)}(\x,\bfl)$,
and
from 12.5 to 5.1 cpu seconds for 1000 evaluations of~$A^{(1)}(\x,\bfl)$.
When the rank deficiency is modest,
the number of multipliers is close to the original number of variables
and we can evaluate about twice as fast.

Setting $\bfl^{(1)} = \bfl$,
$B^{(1)} = B$, and $\bfh^{(1)} = \bfh$,
it is straightforward to formally write
down the $k$th deflated system:
\begin{equation}
   F^{(k)}(\x,\bfl^{(1)},\ldots,\bfl^{(k)}) =
   \left\{
      \begin{array}{lcr}
         F^{(k-1)}(\x,\bfl^{(1)},\ldots,\bfl^{(k-1)}) & = & 0~ \\
         A^{(k-1)}(\x,\bfl^{(1)},\ldots,\bfl^{(k-1)})
        B^{(k)} \bfl^{(k)} & = & 0~ \\
         \multicolumn{1}{r}{\bfh^{(k)} \bfl^{(k)} } & = & 1.
      \end{array}
   \right.
\end{equation}

The benefit of exploiting the structure becomes more significant,
as illustrated by the evaluation of $F^{(2)}$, continuing the
experiment on cyclic 9-roots.  After adding 8 multipliers (corank 2),
in the first deflation, we have 17 variables.  As the corank cannot
decrease, a second deflation will need 16 multipliers.
The time to evaluate $F^{(2)}(\x,\bfl^{(1)},\bfl^{(2)})$
for 1000 random values of $\x$, $\bfl^{(1)}$, and $\bfl^{(2)}$
dropped from 45.1 to 5.66 cpu seconds.

The columns of the Jacobian matrix of $F^{(k)}$ are organized
in $k+1$ blocks.  For example, for $k=2$,
$A^{(2)}(\x,\bfl^{(1)},\bfl^{(2)})$ has in its three block
columns the derivatives with respect to $\x$, $\bfl^{(1)}$,
and~$\bfl^{(2)}$ respectively:
\begin{equation}
  \left[
   \begin{array}{ccc}
      A & \zero & \zero \\
     \left( \frac{\partial A}{\partial \x} \right)
           B^{(1)} \bfl^{(1)} & A B^{(1)} & \zero \\
     \zero & \bfh^{(1)} & \zero \\
     \left( \frac{\partial A^{(1)}}{\partial \x} \right)
           B^{(2)} \bfl^{(2)} &
     \left( \frac{\partial A^{(1)}}{\partial \bfl^{(1)}}
     \right) B^{(2)} \bfl^{(2)} &
     A^{(1)} B^{(2)} \bfl^{(2)} \\
     \zero & \zero & \bfh^{(2)}
   \end{array}
   \right],
\end{equation}
where $A = A(\x)$, $A^{(1)} = A^{(1)}(\x,\bfl^{(1)})$,
and $\frac{\partial}{\partial \bfl^{(1)}}$ means the
same as in~(\ref{eq_derivative_of_matrix}).
Notice that the first three block rows and first two block
columns of $A^{(2)}(\x,\bfl^{(1)},\bfl^{(2)})$ equal
$A^{(1)}(\x,\bfl^{(1)})$.

\section{Applications and Numerical Results}

The implementation has been tested on seventeen examples,
available \newline at
{\tt http://www.math.uic.edu/{\~{}}jan/demo.html}.
The initial approximations for Newton's method were taken from
the end points of solution paths defined by a polynomial homotopy
to find all isolated solutions (see~\cite{Li97,Li03} for a recent surveys).
The numerical results reported in the table below are obtained
with standard machine arithmetic.

Observe the improved numerical conditioning
in Table~\ref{tabresults}.  This observation justifies
the naming\footnote{We are grateful to Erich Kaltofen for
this naming at the ISSAC'04 poster session.}
of our method as a ``re-conditioning'' method.

\begin{table}[hbt]
\begin{center}
\begin{tabular}{|c|c|c||c|c|c|c|} \hline
System & $n$& $m$ & $D$  & $\corank(A(\x))$ & Inverse Condition\# & \#Digits
\\ \hline \hline
 eg1 & 2 & 4 & 1 & 2 $\rightarrow$ 0 & 8.3e-09  $\rightarrow$  5.0e-01 & 9  $\rightarrow$  25\\
 eg2 & 2 & 2 & 1 & 1 $\rightarrow$ 0 & 1.2e-08  $\rightarrow$  1.8e-01  & 9  $\rightarrow$  24\\
 eg3 & 2 & 2 & 1 & 1 $\rightarrow$ 0 & 5.6e-09  $\rightarrow$  1.2e-01  & 9  $\rightarrow$  25\\
 eg4 & 2 & 3 & 2 & 1 $\rightarrow$ 1 $\rightarrow$ 0
 & 3.0e-10 $\rightarrow$  6.4e-02  & 5  $\rightarrow$  15\\
 eg5 & 2 & 4 & 3 & 1 $\rightarrow$ 1 $\rightarrow$ 1 $\rightarrow$ 0
 & 6.4e-11  $\rightarrow$ 5.7e-03  & 4  $\rightarrow$  23\\
\hline
  baker1\cite{KD03} & 2 & 2 & 1 & 1 $\rightarrow$ 0
 & 1.7e-08  $\rightarrow$  3.8e-01  & 9  $\rightarrow$  24\\
\hline
 cbms1\cite{Stu02} & 3 & 11 & 1 & 3 $\rightarrow$ 0
 & 4.2e-05  $\rightarrow$  5.0e-01  & 5  $\rightarrow$  20\\
 cbms2\cite{Stu02} & 3 & 8 & 1 & 3 $\rightarrow$ 0
 & 1.2e-08  $\rightarrow$  5.0e-01  & 8  $\rightarrow$  18\\
 mth191 & 3 & 4 & 1 & 2 $\rightarrow$ 0
 & 1.3e-08  $\rightarrow$  3.5e-02  & 7  $\rightarrow$  13\\
\hline
 decker1\cite{DK80} & 2 & 3 & 2 & 1 $\rightarrow$ 1 $\rightarrow$ 0
 & 3.4e-10  $\rightarrow$  2.6e-02& 6  $\rightarrow$  11\\
 decker2\cite{DKK83}& 2 & 4 & 3 & 1 $\rightarrow$ 1 $\rightarrow$ 1 $\rightarrow$ 0
 & 4.5e-13  $\rightarrow$  6.9e-03 & 5  $\rightarrow$  16\\
 decker3\cite{DKK85}& 2 & 2 & 1 &  1 $\rightarrow$ 0
 & 4.6e-08  $\rightarrow$  2.5e-02 & 8  $\rightarrow$  17\\
\hline
 ojika1\cite{Oji87a} & 2 & 3 & 2 & 1 $\rightarrow$ 1 $\rightarrow$ 0
 & 9.3e-12  $\rightarrow$  4.3e-02 & 5 $\rightarrow$ 12\\
 ojika2\cite{Oji87a} & 3 & 2 & 1 & 1 $\rightarrow$ 0
 & 3.3e-08  $\rightarrow$  7.4e-02 & 6  $\rightarrow$  14\\
 ojika3\cite{OWM83}
 & 3 & 2 & 1 & 1 $\rightarrow$ 0 & 1.7e-08  $\rightarrow$  9.2e-03 & 7  $\rightarrow$  15 \\
 &   & 4 & 1 & 2 $\rightarrow$ 0 & 6.5e-08  $\rightarrow$  8.0e-02 & 6  $\rightarrow$  13\\
 ojika4\cite{Oji88} & 3 & 3 & 2 & 1 $\rightarrow$ 1 $\rightarrow$ 0
 & 1.9e-13  $\rightarrow$  2.4e-04 & 6  $\rightarrow$  11\\
\hline
 cyclic9 & 9 & 4 & 1& 2 $\rightarrow$ 0
 & 5.6e-10  $\rightarrow$  1.8e-03 & 5  $\rightarrow$ 15\\
\hline
\end{tabular}
\caption{Numerical Results, $D$ is the number of deflations needed
         to restore quadratic convergence.  The fifth column
shows the decrease in the corank of the Jacobian matrix for all stages in
the deflation.  In the second to last column we see the change in the estimate
for the inverse condition number of $A(\x)$ at the start of the deflation
to the end of the deflation for~$\x \approx \x^*$.
The last column lists the increase in the number of correct digits from
the initial guess to the final approximation. }
\label{tabresults}
\end{center}
\end{table}

One of the interesting examples
is taken from~\cite{OWM83} and listed as
``ojika3'' in Table~\ref{tabresults}.
This system has two isolated roots: one with multiplicity two,
and the other one has multiplicity four.  Both roots need only
one deflation, but at the double root, the rank of the Jacobian
matrix is two, while the rank is one at the other 4-tuple root.
The program produces two different deflated systems:
one with three multipliers (for the double root) and the other
with two multipliers (for the 4-tuple root).

The so-called cyclic 9-roots problem is the last and largest
example in Table~\ref{tabresults}.
This system is a widely used benchmark in the field of
polynomial system solving, e.g.: \cite{BF91,BF94}, \cite{Fau99},
\cite{GL03}, \cite{GKKTFM04}, \cite{LL01},
with theoretical results in~\cite{Haa96}.
There are 5,594 (333 orbits of size 18) isolated regular cyclic 9-roots,
in addition to the 162 isolated solutions of multiplicity four.
One deflation suffices to restore quadratic convergence on all 162
quadruple roots of this large application.

\section{Conclusion}

Our modified deflation method works in general, is numerically
stable, relatively simple to implement; and perhaps most importantly,
a preliminary implementation on a wide class of examples performs
quite well.

However, the doubling of the number of equations by deflation quickly
leads to huge systems, we will search for ways to limit the number
of deflations.

\bibliographystyle{plain}
\bibliography{newton}

\end{document}